\documentclass[12]{amsart}
\usepackage{amsfonts}
\usepackage{amssymb}
\newtheorem{defi}{Definition}[section]
\newtheorem{satz}[defi]{Theorem}
\newtheorem{lemm}[defi]{Lemma}

\newtheorem{cor}[defi]{Corollary}
\newtheorem{prop}[defi]{Proposition}

\newcommand{\supp}{\operatorname{supp}}

\newcommand{\trace}{\operatorname{trace}}

\newcommand{\h}{\mathcal{H}}
\newcommand{\HH}{\mathcal{H}}

\newcommand{\lspace}{\vspace{0.2cm}}

\newcommand{\DD}{\mathcal{D}}

\newcommand{\C}{\mathbb{C}}
\newcommand{\D}{\mathbb{D}}
\newcommand{\T}{\mathbb{T}}
\newcommand{\N}{\mathbb{N}}

\newcommand{\LL}{\mathcal{L}}
\newcommand{\bmodd}{\operatorname{{\mathrm{BMO^d}}}}
\newcommand{\bmos}{\operatorname{{\mathrm{BMO_{so}}}}}

\newcommand{\bmosd}{\operatorname{{\mathrm{BMO_{so}^d}}}}

\newcommand{\bmond}{\operatorname{{\mathrm{BMO_{norm}^d}}}}
\newcommand{\bmocd}{\operatorname{{\mathrm{BMO^d_{Carl}}}}}

\newcommand{\bmop}{\operatorname{{\mathrm{BMO_{para}}}}}
\newcommand{\bmol}{\operatorname{{\mathrm{BMO_{mult}}}}}
\newcommand{\opf}{\operatorname{{\mathcal{F}_{00}}}}

\newcommand{\FF}{\mathcal{F}}
\newcommand{\ostar}{\circledast}
\newcommand{\honel}{{H^1_\Lambda(\T,S_1)}}
\newcommand{\honemax}{{H^1_{\mathrm{max, d}}(\T,S_1)}}

\newcommand{\mat}{\mathrm{Mat}(\C, n \times n)}

\newcommand{\bmo}{\operatorname{{\mathrm{BMO}}}}
\newcommand{\sbmod}{\operatorname{{\mathrm{SBMO^d}}}}

\begin{document}

\title[Operator BMO spaces]
{Embeddings between operator-valued dyadic BMO spaces}
\author{Oscar Blasco}
\address{Department of Mathematics,
Universitat de Valencia, Burjassot 46100 (Valencia)
 Spain}
\email{oscar.blasco@uv.es}
\author{Sandra Pott}
\address{Department of Mathematics, University of
                    Glasgow, Glasgow G12 8QW, UK}
\email{sp@maths.gla.ac.uk}
\thanks{{\it 2000 Mathematical Subjects Classifications.}
                               Primary 42B30, 42B35, Secondary 47B35 \\
The first author gratefully acknowledges support by the LMS and
Proyectos MTM 2005-08350 and PR2006-0086. The second author
gratefully acknowledges support by EPSRC}
\begin{abstract}
We investigate a scale of dyadic operator-valued BMO spaces,
corresponding to the different yet equivalent characterizations of
dyadic BMO in the scalar case. In the language of operator spaces,
we investigate different operator space structures on the scalar
dyadic BMO space which arise naturally from the different
characterisations of scalar BMO. We also give sharp
dimensional growth estimates for the sweep of functions and its bilinear extension in some of those
different dyadic BMO spaces.
\end{abstract}
\maketitle
\section{Introduction}
Let
$\DD$ denote the collection of dyadic subintervals
of the unit circle $\T$, and let $(h_I)_{I \in \DD}$,
where $h_I = \frac{1}{|I|^{1/2}} ( \chi_{I^+} - \chi_{I^-})$,
be the Haar basis of $L^2(\T)$.
 For $I \in \DD$ and $\phi\in L^2(\T)$, let $\phi_I$ denote the formal Haar
coefficient
$\int_I \phi(t) h_I dt$, and $m_I \phi = \frac{1}{|I|} \int_I \phi(t) dt$
denote the
average of $\phi$ over $I$. We write $P_I(\phi)=\sum_{J\subseteq I}
\phi_Jh_J$.

We say that $\phi\in L^2(\T)$ belongs to dyadic BMO, written
 $\phi\in {\rm BMO^d}(\T)$, if
\begin{equation}\label{bmo1}
   \sup_{I \in \DD}( \frac{1}{|I|} \int_I | \phi(t) - m_I \phi |^2 dt)^{1/2}
                                         < \infty.\end{equation}

Using the identity $P_I(\phi)= (\phi- m_I\phi)\chi_I$, this can
also be written as
\begin{equation}\label{bmo2}
   \sup_{I \in \DD} \frac{1}{|I|^{1/2}} \|P_I( \phi) \|_{L^2}
                                          < \infty,
\end{equation}
or
\begin{equation}  \label{bmo3}
   \sup_{I \in \DD} \frac{1}{|I|} \sum_{J \in \DD, J \subseteq I}
                                                 | \phi_J |^2
                                         < \infty.
\end{equation}
Due to John-Nirenberg's lemma, we have, for $0< p < \infty$, that
$\phi\in {\rm BMO^d}(\T)$ if and only if
\begin{equation}\label{bmo}
   \sup_{I \in \DD}( \frac{1}{|I|} \int_I | \phi(t) - m_I \phi |^p
dt)^{1/p}= \sup_{I \in \DD} \frac{1}{|I|^{1/p}} \|P_I(
\phi)\|_{L^p} <  \infty.\end{equation}

It is well-known that the space  ${\rm BMO^d}(\T)$ has the following equivalent
formulation in terms of boundedness of dyadic paraproducts: The map
\begin{equation}  \label{bmo4}
   \pi_\phi: L^2(\T) \rightarrow L^2(\T), \quad f = \sum_{I \in \DD}
 f_I h_I\mapsto \sum_{I \in \DD}  \phi_I (m_I f) h_I
\end{equation}
defines a bounded linear operator on $L^2(\T)$, if and only if
$\phi\in {\rm BMO^d}(\T)$.

 For real-valued functions, we can also
replace  the boundedness of the dyadic paraproduct
$\pi_\phi$ by the boundedness of its
adjoint operator
\begin{equation}  \label{adjpara}
   \Delta_\phi: L^2(\T) \rightarrow L^2(\T), \quad f = \sum_{I \in \DD}
 f_I h_I\mapsto
                                 \sum_{I \in \DD}  \phi_I f_I
\frac{\chi_I}{|I|}.
\end{equation}
Another equivalent formulation comes from the duality
\begin{equation} \label{h1dual}
{\rm BMO^d}(\T)=(H_d^1(\T))^*,
\end{equation}
 where the dyadic Hardy space $H_d^1(\T)$ consists
of those  functions $\phi\in L^1(\T)$ for which the dyadic square function $\mathcal{S}
\phi =    (\sum_{I\in
\DD}|\phi_I|^2\frac{\chi_I}{|I|})^{1/2}$ is also in $L^1(\T)$.
Let us recall that
  $H_d^1(\T)$ can also be described in terms of dyadic atoms. That is, $H_d^1(\T)$ consists
of functions $\phi= \sum_{k \in
\N} \lambda_k a_k,
              \lambda_k \in \C$,
       $\sum_{k \in \N} | \lambda_k| < \infty $, where the $a_k$
 are  dyadic atoms, i.e. $\supp (a_k)\subset I_k$ for
some $I_k\in \DD$, $\int_{I_k} a_k(t)dt=0$, and
$\|a_k\|_\infty\le\frac{1}{|I_k|}.$
The reader is referred to \cite{M} or to \cite{G} for standard results about
$H^1_d$ and $\bmodd$.

Let
$$S_\phi= (\mathcal{S}\phi)^2=\sum_{I\in
\DD}|\phi_I|^2\frac{\chi_I}{|I|}$$
denote the sweep of the function
$\phi$. Using John-Nirenberg's
lemma, one easily verifies the well-known fact that
\begin{equation} \label{sweep}
\phi\in {\rm BMO^d}(\T) \hbox{ if and only if } S_\phi \in {\rm
BMO^d}(\T).
\end{equation}
The reader is referred to \cite{blasco4} for a proof of
(\ref{sweep}) independent of  John-Nirenberg's lemma.\lspace

The aim of this paper is twofold. Firstly, it is to investigate the spaces of
operator-valued BMO functions  corresponding to characterizations
(\ref{bmo1})-(\ref{h1dual}). In the operator-valued case, these
characterizations are in general no longer equivalent. In the
language of operator spaces, we investigate the different operator
space structures on the scalar space $\bmodd$ which arise
naturally from the different yet equivalent characterisations of
$\bmodd$. The reader is referred to  \cite{BlascoArg,BP,new,psm}
for some recent  results on dyadic BMO and Besov spaces connected
to the ones in this paper. The second aim is to give sharp dimensional estimate
for the operator sweep and its bilinear extension, of which more will be said below,
in these operator $\bmo^d$ norms.

\lspace

We require some further notation for the operator-valued case.
Let $\h$ be a separable, finite or infinite-dimensional Hilbert space.
Let  $\opf$ denote the subspace of
$\LL(\h)$-valued functions
on $\T$
with finite formal Haar expansion.
  Given $e,f\in \h $ and $B \in L^2(\T,\LL(\h))$ we denote by $B_e$ the
function in $L^2(\T,\h)$ defined by $B_e(t)= B(t)(e)$
and by
$B_{e,f}$ the function in $L^2(\T)$ defined by $B_{e,f}(t)=
\langle B(t)(e),f\rangle$.  As in the scalar case,
  let $B_I$ denote the formal Haar coefficients
$\int_I B(t) h_I dt$, and $m_I B = \frac{1}{|I|} \int_I B(t) dt$
denote the average of $B$ over $I$ for any $I \in \DD$. Observe
that for $B_I$ and $m_IB$ to be well-defined operators, we shall
be assuming that the $\LL(\h)$- valued function $B$ is
$weak^*$-integrable. That means, using the duality
$\LL(\h)=(\h\hat\otimes \h)^*$, that $\langle
B(\cdot)(e),f\rangle\in L^1(\T)$ for $e,f\in \h $ and for any
measurable set $A$, there exist $B_A\in \LL(\h)$ such that
$\langle B_A(e),f\rangle=\langle\int_A B(t)(e) dt, f\rangle $ for $e,f\in \h $.

We denote by $\bmodd(\T,\HH)$ the space of Bochner integrable
$\h$-valued  functions $b: \T \rightarrow \h$ such that
\begin{equation}
   \|b\|_{\bmodd}=\sup_{I \in \DD} (\frac{1}{|I|} \int_I \| b(t) - m_I b\|^2
dt)^{1/2}<\infty
\end{equation}
and by $\rm wBMO^d(\T,\h)$ the space of Pettis integrable
$\h$-valued functions $b: \T \rightarrow \h$ such that
\begin{equation}
   \|b\|_{\rm wBMO^d}=\sup_{I \in \DD, e \in \HH, \|e\|=1} (\frac{1}{|I|} \int_I |\langle b(t) - m_I b,
e \rangle|^2 dt)^{1/2}<\infty.
\end{equation}

In the operator-valued case we define the following notions
corresponding to the previous formulations: We denote by
$\bmond(\T,\LL(\h))$ the space of Bochner integrable
$\LL(\h)$-valued  functions $B$ such that
\begin{equation} \label{bmond}
   \|B\|_{\bmond}=\sup_{I \in \DD} (\frac{1}{|I|} \int_I \| B(t) - m_I B
\|^2 dt)^{1/2}<\infty,
\end{equation}
by  ${\rm SBMO^d}(\T,\LL(\h))$ the space of $\LL(\h)$-valued
functions  $B$ such that $B_e\in {\rm BMO^d}(\T,\h)$ for all
$e\in\h$ and
\begin{equation}\label{sbmo}
  \|B\|_{{\rm SBMO^d}}= \sup_{I \in \DD,e \in \HH, \|e \|=1}
                    (\frac{1}{|I|} \int_I \| (B(t) - m_I B)e \|^2
                    dt)^{1/2}< \infty,
\end{equation}
and, finally, by ${\rm WBMO^d}(\T,\LL(\h))$ the space of
$weak^*$-integrable $\LL(\h)$-valued  functions $B$ such that
$B_{e,f}\in {\rm BMO^d}$ for all $e,f\in \h$ and
\begin{multline}\label{wbmo}
\|B\|_{{\rm WBMO^d}}=\sup_{I \in \DD, \|e \|=\|f \|=1}
                    (\frac{1}{|I|} \int_I | \langle(B(t) - m_I B)e,f\rangle
|^2 dt)^{1/2} <\infty,
\end{multline}
or, equivalently, such that
$$
\|B\|_{{\rm WBMO^d}}= \sup_{e \in \HH, \|e \|=1} \|B_e\|_{\rm
wBMO^d(\T,\HH)} =\sup_{A \in S_1, \| A \|_1 \le 1 } \| \langle B,
A \rangle \|_{\bmodd(\T)} < \infty.
$$
Here, $S_1$ denotes the ideal of trace class operators in $\LL(\h)$,
and $\langle B, A \rangle$ stands for the scalar-valued function
given by $\langle B, A \rangle (t) = \trace( B(t) A^*)$.

The space $ \bmocd (\T,\LL(\h))$ is the space of
$weak^*$-integrable operator-valued functions for which
\begin{equation} \label{def:bmocd}
   \|B\|_{\bmocd}=\sup_{I \in \DD} (\frac{1}{|I|}
\sum_{J \in \DD, J \subseteq I}  \| B_J \|^2 )^{1/2} < \infty.
\end{equation}

We would like to point out that while $B$ belongs to one of the
spaces $ \bmond(\T,\LL(\h)),{\rm WBMO^d}(\T,\LL(\h))$) or
$B\in\bmocd(\T,\LL(\h))$  if and only if $B^*$ does, this is  not
the case for the space $\mathrm{SBMO}^d(\T,\LL(\h))$. This leads
to the following notion
 (see \cite{gptv2, petermichl, pxu}):  We say that
 $B\in\bmosd(\T,\LL(\h))$, \label{bmosd} if
 $B$ and $B^*$ belong to $ {\rm SBMO^d}(\T,\LL(\h))$.
We define \begin{equation} \label{def:bmoso}\|B\|_{\bmosd}=
\|B\|_{{\rm SBMO^d}}+\|B^*\|_{{\rm SBMO^d}}.\end{equation}

We now define another operator-valued BMO space, using the notion of Haar
multipliers.

 A sequence $(\Phi_I)_{I \in \DD}$, $\Phi_I\in L^2(I,\LL(\h))$
for all $I\in \DD$, is said to be an \emph{operator-valued Haar
multiplier} (see \cite{per, BP}), if there exists $C>0$ such that
$$\|\sum_{I\in \DD}\Phi_I(f_I)h_I\|_{L^2(\T,\h)}\le C (\sum_{I\in
\DD}\|f_I\|^2)^{1/2} \text{ for all }
(f_I)_{I \in \DD} \in l^2(\DD,\h).$$
We write  $\|(\Phi_I)\|_{mult}$ for the norm of the corresponding operator
on
$L^2(\T,\h)$.

Letting  again
as in the scalar-valued case
$P_I B =\sum_{J\subseteq I} h_JB_J$,
we denote the space of those $weak^*$-integrable
$\LL(\h)$-valued functions for
which $(P_IB)_{I\in\DD}$ defines a bounded operator-valued Haar multiplier on $L^2(\T, \h)$
by
   $ \bmol (\T,\LL(\h))$ and write
 \begin{equation}\label{bmol}\|B\|_{\bmol}=
\|(P_IB)_{I\in\DD}\|_{mult}.
\end{equation}
We shall use the notation $\Lambda_B(f)=\sum_{I \in \DD} (P_I B) (f_I) h_I$.

Let us mention that there is a further BMO space, defined in terms
of paraproducts, which is very much connected with
$\bmol(\T,\LL(\h))$ and was studied in detail in \cite{new}.
Operator-valued paraproducts are of particular interest, because
they can be seen as dyadic versions of vector Hankel operators or
of vector Carleson embeddings, which are important in the real and
complex analysis of matrix valued functions and its applications in the theory
of infinite-dimensional linear systems (see e.g.~\cite{jp}, \cite{jpp1}).

Let $B \in \opf$. We define the dyadic operator-valued paraproduct
with symbol $B$,
$$
   \pi_B: L^2(\T, \h) \rightarrow L^2(\T, \h), \quad f = \sum_{I \in \DD}
 f_I h_I\mapsto
                                 \sum_{I \in \DD}  B_I (m_I f) h_I,
$$
and
$$
   \Delta_B: L^2(\T, \h) \rightarrow L^2(\T, \h), \quad f = \sum_{I \in \DD}
 f_I h_I\mapsto
                                 \sum_{I \in \DD}  B_I (f_I)
\frac{\chi_I}{|I|}.
$$
It is easily seen that $(\pi_B)^* = \Delta_{B^*}$.

We denote  by $ \bmop (\T,\LL(\h))$ the space of
$weak^*$-integrable operator-valued functions for which
$\|\pi_{B}\| < \infty$
 and write \begin{equation}
\label{bmop}\|B\|_{\bmop}= \|\pi_B\|.
\end{equation}
  We refer the reader to \cite{BlascoArg, new} and \cite{mei1, mei}
 for results on this space.
It is  elementary to see that
\begin{equation}  \label{mult}
  \Lambda_B( f )=
                                 \sum_{I \in \DD}  B_I (m_I f) h_I
                            + \sum_{I \in \DD} B_I (f_I) \frac{\chi_I}{|I|} = \pi_B f + \Delta_B f.
\end{equation}
Hence $\Lambda_B=\pi_B+\Delta_B$ and
$(\Lambda_B)^*=\Lambda_{B^*}$. This shows that
$\|B\|_{\bmol}=\|B^*\|_{\bmol}$.

 Let us finally
denote by ${\rm BMO_{spara}}(\T,\LL(\h))$ the space of symbols $B$
such that $\pi_{B}$ and $\pi_{B^*}$ are bounded operators, and
define
\begin{equation}
\label{bmosp}\|B\|_{\rm BMO_{spara}}= \|\pi_B\|+\|\pi_{B^*}\|.
\end{equation}
Since   $\Delta_B=\pi^*_{B^*}$, one concludes  that ${\rm
BMO_{spara}}(\T,\LL(\h))\subseteq \bmol
 (\T,\LL(\h))$.

We write $\approx$ for equivalence of norms up to a constant (independent of the
dimension of the Hilbert space $\h$, if this appears), and similarly $\lesssim, \gtrsim$ for
the corresponding one-sided estimates up to a constant.

Recall that for a given Banach space $(X, \| \cdot \|)$, a family
of norms $( M_n(X), \| \cdot\|_n)$ on the spaces $M_n(X)$ of
$X$-valued $n \times n$ matrices defines an \emph{operator space
structure} on $X$, if $\|\cdot\|_1 \approx \|\cdot\|$,
\begin{enumerate}
\item[(M1)] $\| A \oplus B \|_{n +m} = \max \{ \|A\|_n, \|B\|_m \}$ for $A \in M_n(X)$,
          $B \in M_m(X)$
\item[(M2)] $ \|\alpha A \beta \|_{m} \le \|\alpha \|_{M_{n,m}(\C)}
\| A \|_n \|\beta\|_{M_{m,n}(\C)} $ for all $A \in M_n(X)$
 and all scalar matrices  $\alpha \in M_{n,m}(\C)$, $\beta \in M_{m,n}(\C)$.
\end{enumerate}
(see e.~g.~\cite{effros}). One verifies easily that all the $\bmodd$-norms
on $\LL(\h)$-valued functions defined above,
except $\bmond$ and $\bmocd$, define operator space
structures on $\bmodd(\T)$ when taken for n-dimensional $\h$, $n \in \N$. \lspace

The aim of the paper  is to show the following strict inclusions for infinite-dimensional
$\h$:
\begin{multline}   \label{eq:inclchain}
 \bmond(\T,\LL(\h))
   \subsetneq \bmol(\T,\LL(\h))\subsetneq\\ \subsetneq {\bmosd}(\T,\LL(\h))
   \subsetneq
{\rm WBMO^d}(\T,\LL(\h))
\end{multline}
and
\begin{equation} \label{carles}
\bmocd(\T,\LL(\h))  \subsetneq {\rm
BMO_{spara}}(\T,\LL(\h))\subsetneq
\bmol(\T,\LL(\h)).\end{equation}
 This means that the
corresponding inclusions of operator spaces over $\bmodd(\T)$, where they apply, are
completely bounded, but not completely isomorphic (for the
notation, see again e.~g.~\cite{effros}). We will also consider the preduals for some of the spaces
shown. Finally, we will give sharp estimates for the dimensional growth of the sweep and its bilinear
extension on $\bmop$, $\bmol$ and $\bmond$, completing results in \cite{new} and \cite{mei}.

The paper is organized as follows. In Section 2, we prove the chains
of strict inclusions (\ref{eq:inclchain}) and (\ref{carles}).
Actually the only nontrivial inclusion to be shown is $\bmond(\T,\LL(\h))
   \subset \bmol(\T,\LL(\h))$.
For this purpose, we introduce a new Hardy space $H^1_\Lambda$ adapted to the problem, and then the result
 can be shown from an estimate on the dual side. The remaining inclusions are
immediate consequences of the definition, and only the
counterexamples showing that none of the spaces are equal need to be found.

The reader is referred to \cite{mei1} for more on the theory of
operator-valued Hardy spaces.

Section 3 deals with dimensional growth properties of the \emph{operator sweep}
and its bilinear extension. We define the operator sweep for $B \in \opf$,
$$
    S_B = \sum_{I \in \DD} \frac{\chi_I}{|I|} B_I^* B_I,
$$
and its bilinear extension
$$
    \Delta[U^*,V]= \sum_{I \in \DD} \frac{\chi_I}{|I|} U_I^* V_I   \qquad (U,V \in \opf).
$$
 These maps are of interest for several reasons. They are closely connected with
the paraproduct and certain bilinear paraproducts,
they provide a tool to understand the dimensional growth in the John-Nirenberg lemma,
and they are useful to understand
products of paraproducts and products of certain other operators (see \cite{new}, \cite{psm}).

Considering (\ref{sweep}) in the operator valued case, it was shown
in \cite{new} that
\begin{equation}\label{normsweep} \|S_B\|_{\rm
BMO^d_{mult}}+\|B\|^2_{\rm SBMO^d}\approx \|B\|^2_{\rm
BMO^d_{para}}.
\end{equation}
Here, we prove the bilinear analogue
\begin{equation}\label{normbisweep} \|\Delta[U^*,V]\|_{\rm
BMO^d_{mult}}+\sup_{I \in \DD} \frac{1}{|I|} \|\sum_{J \subset I} U_J^* V_J\|
     \approx \|\pi_U^* \pi_V\|.
\end{equation}
It was also shown in \cite{new} that
\begin{equation}\label{estsweep}
\|S_B\|_{\rm SBMO^d}\le C \log(n+1) \|B\|^2_{\rm SBMO^d}
\end{equation}
for $dim(\HH)=n$, where $C$ is a constant independent of $n$, and that this estimate is sharp.

We extend this by proving sharp estimates
of $\|S_B\|$ and $\|\Delta[U^*,V]\|$ in terms of $\|B\|, \|U\|, \|V\|$ with respect to the
norms in ${\rm SBMO}^d$, $\bmop$, $\bmol$ and $\bmond$.

\section{Strict inclusions}

Let us start by stating the following characterizations of
$\mathrm{SBMO}$ to be used later on. Some of the
equivalences  can be found in \cite{gptv2}, we give the proof for
the convenience of the reader.

\begin{prop}\label{carbmoso} Let $B\in{\rm
SBMO^d}(\T,\LL(\h))$. Then

\begin{eqnarray*}\|B\|^2_{{\rm SBMO^d}} &=& \sup_{e \in \HH,
\|e
\|=1}
\|B_e\|^2_{\bmodd(\T,\HH) }\\
&=&\displaystyle\sup_{I\in \DD,\|e\|=1}\frac{1}{|I|}\|P_I(
B_e)\|^2_{L^2(\h)} \\
&=&\displaystyle\sup_{I\in \DD}\frac{1}{|I|}\|\sum_{J\subseteq I}
B_J^* B_J\| \\\
&=&\displaystyle
   \sup_{I \in \DD}\left\|
    \frac{1}{|I|} \int_I (B(t) - m_I B)^* (B(t) - m_I B) dt
\right\|\\
&=&\sup_{I\in \DD} \|m_I(B^*B)-m_I(B^*)m_I(B)\|.
\end{eqnarray*}
\end{prop}
\proof The two first equalities are obvious from the definition.
Now observe  that
$$\|\sum_{J\subseteq I} B_J^*
B_J\|=\sup_{\|e\|=1,
\|f\|=1}\sum_{J\subseteq I} \langle B_J(e),
B_J(f)\rangle=\sup_{\|e\|=1}\sum_{J\subseteq I} \| B_J(e)
\|^2=\|P_I(B_e)\|^2_{L^2(\h)}.$$

The other equalities follow from
\begin{eqnarray*}
  \|m_I(B^*B)-m_I(B^*)m_I(B)\| &=&  \left\|\frac{1}{|I|} \int_I
(B(t)-m_I B)^* (B(t) - m_IB) dt  \right\| \\
  &=& \sup_{e \in \HH, \|e\|=1 }\frac{1}{|I|} \int_I  \langle
(B(t)-m_I B)^* (B(t) - m_IB)e,e \rangle dt \\
&=&
  \sup_{e \in \HH, \|e\|=1 }\frac{1}{|I|} \int_I \| P_I B e\|^2 dt.
\end{eqnarray*}
\qed

\lspace

\begin{lemm} Let $B=\sum_{k=1}^N B_k r_k$ where $ r_k=\sum_{|I|=2^{-k}}|I|^{1/2} h_I$ denote the Rademacher
functions.  Then
\begin{equation}\label{n1}
\|B\|_{\sbmod}=
\sup_{\|e\|=1}(\sum_{k=1}^N \|B_k e\|^2)^{1/2}
\end{equation}
\begin{equation}\label{n2}
\|B\|_{\bmos}=
\sup_{\|e\|=1}(\sum_{k=1}^N \|B_k
e\|^2)^{1/2}+\sup_{\|e\|=1}(\sum_{k=1}^N\|B^*_k e\|^2)^{1/2}
\end{equation}
\begin{equation}\label{n3}
\|B\|_{\mathrm{WBMO}^d}= \sup_{\|f\|=\|e\|=1}
(\sum_{k=1}^N |\langle B_k e,f\rangle|^2
)^{1/2}.
\end{equation}
\end{lemm}
\proof This follows from standard Littlewood-Paley theory.
\qed

\lspace

For $x,y\in \h$ we denote by $x\otimes y$ the rank 1 operator in
$\LL(\h)$ given by $(x\otimes y)(h)=\langle h,y\rangle x$. Clearly
$(x\otimes y)^*=(y\otimes x)$.

\begin{prop} \label{firstinc} Let $\dim \h = \infty$. Then
$$\bmol \subsetneq\bmosd(\T, \LL(\HH)) \subsetneq \sbmod(\T,\LL(\h))\subsetneq {\rm
WBMO^d}(\T,\LL(\h)).$$
\end{prop}
\proof Note that if $(\Phi_I)_{I\in\DD}$ is a Haar multiplier then
\begin{equation}\label{debilmult}
\sup_{I\in \DD, \|e\|=1}  |I|^{-1/2}
\|\Phi_I(e)\|_{L^2(\T,\h)}\le \|(\Phi_I)\|_{mult}.
 \end{equation}

The first inclusion thus follows from (\ref{debilmult}) and Proposition
\ref{carbmoso}. The other inclusions are immediate.
Let us see that they are strict. It was shown in \cite{gptv2} that $\bmol(\T,\LL(\h)) \neq
\bmosd(\T,\LL(\h))$.

 Let $(e_k)$ is an orthonormal basis of $\h$ and
$h\in \h$ with $\|h\|=1$. Hence by (\ref{n1}),
$B=\sum_{k=1}^\infty h\otimes e_k
\ r_k\in \sbmod $ and $B^*=\sum_{k=1}^\infty e_k\otimes h \
r_k\notin\sbmod(\T,\LL(\h))$.
Thus $B\in\sbmod(\T, \LL(\HH)) \setminus \bmosd(\T,\LL(\h))$. Similarly by (\ref{n1}) and (\ref{n3}), $B\in
 {\rm
WBMO^d}(\T,\LL(\h))\setminus \sbmod(\T,\LL(\h))$. \qed

\lspace

Note that
\begin{equation}\label{form}
\Lambda_B f= B  f -\sum_{I\in \DD}(m_IB)(f_I) h_I
\end{equation}
which allows to conclude immediately that
 $L^\infty(\T,\LL(\h))\subseteq \bmol(\T,\LL(\h))$.

Our next objective is to see that $\bmond(\T,\LL(\h)) \subsetneq
\bmol(\T,\LL(\h))$. For that,
 we need again some more notation.

Let $S_1$ denote the ideal of trace class operators on $\h$ and
recall that $S_1=\h\hat\otimes\h$ and $(S_1)^*=\LL(\h)$ with the
pairing $\langle U,(e \otimes d)\rangle= \langle U(e),  d\rangle.$

  It is easy to see that
the space $\bmol(\T,\LL(\h))$ can be embedded isometrically into the dual of
a certain
$H^1$ space
of $S_1$ valued functions:
\begin{defi}  Let $f,g\in L^2(\T,\h)$.  Define$$
   f \ostar g =  \sum_{I \in \DD} h_I (f_I \otimes m_I g + m_I f
\otimes g_I).
$$
Let $\honel$  be the space of functions $ f=\sum_{k=1}^\infty \lambda_k f_k
\ostar g_k$ such that $f_k,g_k\in L^2(\T,\h)$,
$\|f_k\|_{2}=\|g_k\|_2=1$ for all $k \in \N$, and $\sum_{k=1}^\infty
|\lambda_k| <\infty.$

We endow the space with the norm given by the infimum of $\sum_{k=1}^\infty
|\lambda_k|$ for all possible decompositions.

\end{defi}

\lspace

With this notation, $B \in \bmol$ acts on $f \ostar g$ by
$$
\langle B,  f \ostar g \rangle
= \int_\T \langle B(t),  (f  \ostar g)(t) \rangle dt =
\langle \Lambda_B f,g \rangle.
$$
By definition of $\honel$, $\| B\|_{(\honel)^*} = \|\Lambda_B\|$.

We will now define a further  $H^1$ space of $S_1$-valued functions.
For $F \in L^1(\T, S_1)$, define the dyadic Hardy-Littlewood maximal
function $F^*$
of $F$ in
the usual way,
$$
     F^*(t) = \sup_{I \in \D, t \in I} \frac{1}{|I|} \int_I \| F(s)\|_{S_1}
ds.
$$
Then let $\honemax$ be given by functions $ F \in L^1(\T,S_1)$
such that $ F^* \in L^1(\T) $.  By a result of Bourgain
(\cite{bourgain}, Th.12), $\bmond$ embeds continuously into
$(\honemax)^*$ (see also \cite{blasco1,blasco3}).

\begin{lemm} \label{dual}
   $\honel \subseteq \honemax.$
\end{lemm}
\proof It is sufficient to show that there is a constant $C >0$ such
that for all $f,g \in L^2(\T,\h)$, $f \ostar g \in \honemax$, and
$\| f \ostar g \|_{\honemax} \le C \|f\|_2 \|g\|_2$.
One verifies that
$$
   f \ostar g = \sum_{I \in \DD} h_I (f_I \otimes m_I g + m_I f
\otimes g_I)= f \otimes g - \sum_{I \in \DD}\frac{\chi_I}{|I|} f_I \otimes
g_I.
$$
Towards the estimate of the maximal function,
let $E_k$ denote the expectation with respect to the
$\sigma$-algebra generated
by dyadic intervals of length $2^{-k}$,
$$
    E_k F = \sum_{I \in \DD, |I| > 2^{-k}} h_I F_I,
$$
for each $k \in \N$.
Then we have
\begin{equation}
   E_k( f \ostar g) = (E_k f) \ostar (E_k g),
\end{equation}
as
$$
     \sum_{I \in \DD, |I| > 2^{-k}} h_I (f_I \otimes m_I g + m_I f \otimes
g_I)
    = \sum_{I \in \DD} h_I ((E_k f)_I \otimes m_I (E_k g) +
       m_I (E_k f) \otimes (E_k g)_I).
$$
Thus
\begin{multline*}
    (f \ostar g)^*(t) = \sup_{k \in \N} \|E_k(f \ostar g)(t)\|_{S_1}
      \le \sup_{k \in \N} \|(E_k f)(t)\| \|(E_k g)(t)\| +\sum_{I \in \DD}
\frac{\chi_I(t)}{|I|} \|f_I\| \|g_I\|\\
      \le \|f^*(t)\| \|g^*(t)\| +\sum_{I \in \DD} \frac{\chi_I(t)}{|I|}
\|f_I\| \|g_I\|,
\end{multline*}
and
$$
    \|(f \ostar g)^*\|_1 \le \|f^*\|_2 \|g^*\|_2 + \|f\|_2 \|g\|_2 \le C
\|f\|_2 \|g\|_2
$$
by the Cauchy-Schwarz inequality and boundedness of the dyadic
Hardy-Littlewood maximal function on $L^2(\T, \h)$.
\qed

\lspace

\noindent
In particular, $\honel \subseteq L^1(\T,S_1)$.

\lspace

We can now prove our inclusion result:
\begin{satz}\label{maininclu}
    $\bmond(\T,\LL(\h)) \subsetneq \bmol(\T,\LL(\h))$.
\end{satz}

\proof The inclusion follows by Lemma \ref{dual}, duality and Bourgain's
result.

To see that the spaces do not coincide, use the fact that
$\bmodd(\ell_\infty)\subsetneq \ell_\infty(\bmodd)$ to find for
each $N \in \N$ functions $b_k \in \bmo$, $k=1,...,N$, such that
$\sup_{1\le k\le N}\|b_k\|_{\bmodd}\le 1$, but
$\|(b_k)_{k=1,\dots,N}\|_{\bmo^d(\T,l^\infty_N)}\ge c_N$, $c_N
{\rightarrow} \infty$ as $N \to \infty$.

 Let $(e_k)_{k \in \N}$ be
an orthonormal basis of $\HH$, and
consider the operator-valued function
$B(t)=\sum_{k=1}^{N} b_k(t)e_k\otimes e_k\in L^2(\T,\LL(\ell_2))$. Clearly
$B_I= \sum_{k=1}^N (b_k)_I e_k \otimes e_k$, and
for each $\C^N$-valued function $f= \sum_{k=1}^N f_k e_k$, $f_1, \dots, f_N \in L^2(\T)$, we have
$$\Lambda_B(f)=\sum_{k=1}^N \Lambda_{b_k}(f_k)e_k . $$
Choosing the $f_k$ such that $\|f\|_2^2=\sum_{k=1}^N
\|f_k\|^2_{L^2(\T)}=1$, we find that
$$
\|\Lambda_B(f)\|^2_{L^2(\T,\ell_2)}=\sum_{k=1}^N
\|\Lambda_{b_k}(f_k)\|^2_{L^2(\T)}\le C \sum_{k=1}^N
\|{b_k}\|^2_{\bmodd}\|f_k\|^2_{L^2(\T)}\le C,
$$
where $C$ is a constant independent of $N$. Therefore, $\Lambda_B$ is bounded.

But since
$\|B\|_{\bmond}=\|(b_k)_{k=1,\dots,N}\|_{{\bmodd}(\T,l^\infty_N)}\ge
c_N$, it follows that  $\bmol(\T)$ is not continuously embedded in
$\bmond(\T,
\LL(\h))$.
From the open mapping theorem, we obtain inequality of the spaces.
\qed
\medskip

The next proposition shows that the space $\bmocd$ belongs to a
different scale than $\bmond$ and $\bmol$.
\begin{prop}
$L^\infty(\T,\LL(\h))  \nsubseteq \bmocd(\T,\LL(\h)).$
\label{subs:inf-carl}
\end{prop}
\proof This follows from the result $L^\infty(\T,
\LL(\HH))\nsubseteq \bmop $ in \cite{mei} (see Lemma \ref{mei}
below) and next proposition. We give a simple direct argument.
Choose an orthonormal basis of $\h$ indexed by the elements of
$\DD$, say $(e_I)_{I \in \DD}$, and let $\Phi_I = e_I \otimes
e_I$, $\Phi_I h = \langle h, e_I \rangle e_I$. Let $\lambda_I =
|I|^{1/2}$ for $I \in \DD$, and define $B = \sum_{I \in \DD} h_I
\lambda_I \Phi_I$. Then $\sum_{I \in \DD} \|B_I\|^2 = \sum_{I \in
\DD} |I| = \infty$, so in particular $B \notin
\bmocd(\T,\LL(\h))$. But the operator function $B$ is diagonal
with uniformly bounded diagonal entry functions $\phi_I(t)
=\langle B(t) e_I, e_I \rangle = |I|^{1/2} h_I(t)$, so $B \in
L^\infty(\LL(\h))$.\qed

\begin{prop} \label{subs:carl-para}
$$\bmocd(\T,\LL(\h))  \subsetneq {\rm
BMO_{spara}}(\T,\LL(\h))\subsetneq {\rm BMO_{mult}}(\T,\LL(\h)).$$

\end{prop}
\proof The inclusion $\bmocd \subseteq {\rm BMO_{spara}}$ is easy,
 since (\ref{def:bmocd}) implies that for $B \in \bmocd$,
the $\bmocd$ norm equals the norm of the scalar $\bmodd$ function
given by $|B|:=\sum_{I \in \DD} h_I \|B_I\|$. For $f \in L^2(\h)$,
let $|f|$ denote the function given by $|f|(t) = \|f(t)\|$. Thus
$$
  \|\pi_B f \|_2^2 = \sum_{I \in \DD} \|B_I m_I f\|^2 \le
                     \sum_{I \in \DD} (\|B_I\| m_I |f|)^2 = \| \pi_{|B|} |
f|\|.
$$
The boundedness of $\pi_{B^*}$ follows analogously.

To show that $\bmocd \neq {\rm BMO_{spara}}$, we can use the
diagonal operator function $B$ constructed in Proposition
\ref{subs:inf-carl}. There, it is shown that $B \notin \bmocd$,
and that the diagonal entry functions $\phi_I = \langle Be_I, e_I
\rangle$ are uniformly bounded. Since the paraproduct of each
scalar-valued $L^\infty$ function is bounded, we see that $\pi_B =
\bigoplus_{I \in \DD} \pi_{\phi_I}$ is bounded. Similarly,
$\pi_{B^*}$ is bounded. Thus $B \in {\rm BMO_{spara}}$.
It is clear from (\ref{mult}) that ${\rm BMO_{spara}}(\T,\LL(\h))\subseteq {\rm
BMO_{mult}}(\T,\LL(\h))$.

 Using that
$L^\infty(\T,\LL(\h))\nsubseteq {\rm BMO_{spara}}(\T,\LL(\h))$
(see \cite{mei}), one concludes that ${\rm
BMO_{spara}}(\T,\LL(\h))\neq {\rm BMO_{mult}}(\T,\LL(\h))$. \qed
\medskip

\section{Sharp dimensional growth of the sweep}
We begin with the following lower estimate of the $\bmop$ norm in terms of the $L^\infty$ norm of
certain $\mat$-valued functions from \cite{mei}.
\begin{lemm}\label{mei} (see \cite{mei}, Thm 1.1.) There exists an absolute constant
$c >0$ such that for each $n\in \N$, there exists a measurable
 function $F:\T \rightarrow \mat$ with $\|F\|_\infty \le 1$ and
$\|\pi_F\| \ge c \log(n+1)$.
\end{lemm}

Here are our dimensional estimates of the sweep.
\begin{satz} There exists an absolute constant $C >0$ such that for
each $n \in \N$ and each measurable function $B: \T \rightarrow  \mat$,
\begin{equation}  \label{eq:sharppara}
    \| S_B \|_{\bmop} \le C \log(n+1) \| B\|^2_{\bmop},
\end{equation}
\begin{equation}  \label{eq:sharpmult}
    \| S_B \|_{\bmol} \le C (\log(n+1))^2 \| B\|^2_{\bmol},
\end{equation}
\begin{equation}  \label{eq:sharpnorm}
    \| S_B \|_{\bmond} \le C (\log(n+1))^2 \| B\|^2_{\bmond},
\end{equation}
and the dimensional estimates are sharp.
\end{satz}
\proof Let $B: \T \rightarrow  \mat$ be measurable. Since
$\|S_B\|_* = \lim_{k \to \infty} \|S(E_k B) \|_*$ in all of the
above BMO norms ( because we are in the finite-dimensional
situation) it suffices to consider the case $B \in \FF_{00}$.

We start by proving (\ref{eq:sharppara}). Since
\begin{equation}   \label{eq:paraest}
\|\pi_B\| \le C'
\log(n+1) \|B\|_{\bmosd}
\end{equation}
 for some absolute constant $C' >0$
(see \cite{ntv}, \cite{katz}) and
\begin{equation}  \label{eq:strongmult}
   \|B\|_{\bmosd} \le \|B\|_{\bmol},
\end{equation}
we have
\begin{equation*}  \|S_B \|_{\bmop} \le C'
\log(n+1) \|S_B\|_{\bmol} \le C \log(n+1) \|B\|^2_{\bmop}
\end{equation*}
by (\ref{normsweep}).

For the sharpness of the estimate, take $F$ as in Lemma \ref{mei}.  Again,
approximating by $E_k F$, we can assume that $F \in \FF_{00}$.
Since each function in $L^\infty(\T, \mat)$ is the linear
combination of 4 nonnegative-matrix valued functions, the $L^\infty$-norm of which is controlled by the
norm of the original function, we can (by replacing $c$ with a smaller constant) assume
that $F$ is a nonnegative matrix-valued function. Each such nonnegative
matrix-valued function $F$ can be written as $F = S_B$ with $B \in
\FF_{00}$, for example by choosing $B = \sum_{I \in \DD, |I|=
{2^{-k}}} h_I B_I$, where $B_I = |I|^{1/2} (F^I)^{1/2}$, $F =
\sum_{I \in \DD, |I|= {2^{-k}}} \chi_I F^I$. It follows that
\begin{multline*}
 \|S_B \|_{\bmop}\ge c \log(n+1) \|S_B\|_\infty  \\
   \ge c/2 \log(n+1)( \|S_B\|_{\bmol} + \|B\|^2_{\bmosd})
      \gtrsim \log(n+1) \|B\|^2_{\bmop}
\end{multline*}
again by (\ref{normsweep}). Here, we use the estimate
$\|B\|^2_{\bmosd} \le \|S_B\|_\infty$, which can easily be obtained by
$$
    \|P_I B e\|_2^2 = \|S_{P_I Be}\|_1 \le |I| \|S_{P_I Be}\|_\infty \le |I| \|S_{P_I B}\|_\infty
   \le  |I| \|S_{ B}\|_\infty \text{ for } e \in \h, \|e\|=1.
$$
This proves that (\ref{eq:sharppara}) is sharp.

Let us now show (\ref{eq:sharpmult}). Note that by
(\ref{normsweep}) and (\ref{eq:paraest}), for $B \in \FF_{00}$,
$$
     \|S_B \|_{\bmol} \lesssim \|B\|^2_{\bmop} \le {C}'^2 \log(n+1)^2 \|B\|^2_{\bmol}.
$$

For sharpness, choose $B \in \FF_{00}$, $\|B \|_{\infty} \le 1$,
$\|\pi_B\|\ge c \log(n+1)$  as  above, to obtain
\begin{multline*}
    \|S_B \|_{\bmol} + \|B\|^2_{\bmosd}  \gtrsim \|B\|^2_{\bmop}  \\
\ge c^2 \log(n+1)^2 \|B\|_\infty^2
                            \ge  c^2 \log(n+1)^2 \|B\|_{\bmol}^2
\end{multline*}
and thus
$$
\|S_B \|_{\bmol} \gtrsim  \log(n+1)^2 \|B\|_{\bmol}^2,
$$
as $\|B\|_{\bmosd} \le \|B\|_{\bmol}$.

Finally, let us show (\ref{eq:sharpnorm}). Again, we can restrict
ourselves to the case $B \in \opf$ by an approximation argument.
We use the fact that the UMD constant of $\mat$ is equivalent to
$\log(n+1)$ (see for instance \cite{Pi1}) and the representation
$$
   S_B(t) = \int_{\Sigma} (T_\sigma B)^*(t) (T_\sigma B)(t)  d \sigma \qquad(B \in \opf)
$$
(see \cite{new}, \cite{gptv2}),
where $T_\sigma$ denotes the dyadic martingale transform $B \mapsto T_\sigma B = \sum_{I \in \DD} \sigma_I h_I B_I$,
$\sigma= (\sigma_I)_{I \in \DD} \in\{-1,1\}^\DD$,
and $d \sigma$ the natural product probability  measure on $\Sigma =\{-1,1\}^\DD$ assigning measure $2^{-n}$ to
cylinder sets of length $n$,
to prove that
\begin{multline*}
   \|P_I S_B\|_{L^1(\T, \mat)} = \|P_I S_{P_I B}\|_{L^1(\T, \mat)}
     \le 2 \|S_{P_I B}\|_{L^1(\T, \mat)}  \\
  \lesssim (\log(n+1))^2 \|P_I B\|_{L^2(\T, \mat)}^2 \le (\log(n+1))^2
                         |I| \|B\|^2_{\bmond},
\end{multline*}
which gives the desired inequality.

To prove sharpness, choose $B \in \FF_{00}$, $\|B \|_{\infty} \le 1$, $\|\pi_B\|\ge c \log(n+1)$
and note that by Theorem \ref{maininclu},
\begin{multline*}
    \|S_B\|_{\bmond} + \|B \|^2_{\bmosd} \gtrsim \|S_B\|_{\bmol} + \|B \|^2_{\bmosd} \\
            \gtrsim \|B\|^2_{\bmop} \ge c^2 \log(n+1)^2 \|B\|^2_\infty \ge c^2  \log(n+1)^2 \|B\|^2_{\bmond}.
\end{multline*}
Since $\|B \|_{\bmosd} \le \|B\|_{\bmond}$, this implies
$$
\|S_B\|_{\bmond}  \gtrsim \log(n+1)^2 \|B\|^2_{\bmond}.
$$
\qed

We now consider the bilinear extension of the sweep. By
\cite{psm}, \cite{new} or \cite{BlascoArg}
\begin{equation} \label{eq:bidentity}
\pi_{U}^* \pi_{V} = \Lambda_{\Delta[U^*,V]} + D_{U^*,V} \qquad(U,V \in \opf),
\end{equation}
where $D_{U^*,V}$ is given by $D_{U^*,V} h_I e = h_I
\frac{1}{|I|}\sum_{J \subset I} U^*_J V_J e$ for $I \in \DD$, $e
\in \h$.
\begin{prop} \label{prop:bisweep}
$$
      \|\pi_{U}^* \pi_{V}\| \approx \|\Delta[U^*,V]\|_{\bmol} + \sup_{I \in \DD}\frac{1}{|I|}\|\sum_{J \subset I} U^*_J V_J \|
\qquad (U,V \in \opf).
$$
\end{prop}
\proof Obviously $\|D_{U^*,V}\| = \sup_{I \in
\DD}\frac{1}{|I|}\|\sum_{J \subset I} U^*_J V_J \|$. Thus by
(\ref{eq:bidentity}),
$$
\|\pi_{U}^* \pi_{V}\| \le \|\Delta[U^*,V]\|_{\bmol} + \sup_{I \in
\DD}\frac{1}{|I|}\|\sum_{J \subset I} U^*_J V_J \|.
$$
For the reverse estimate, it suffices to observe that $D_{U^*,V}$ is the block diagonal of the
operator $\pi_{U}^* \pi_{V}$ with respect to the orthogonal subspaces $h_I \h$, $I \in \DD$
and therefore $\|D_{U^*,V}\| \le \|\pi_{U}^* \pi_{V} \|$.
\qed

\lspace
\noindent
Here are the dimensional estimates of the bilinear map $\Delta$.
\begin{cor} There exists an absolute constant $C >0$ such that for
each $n \in \N$ and each pair of measurable functions $U,V: \T \rightarrow  \mat$,
\begin{equation}  \label{eq:bistrong}
    \| \Delta[U^*,V] \|_{\sbmod} \le C \log(n+1) \|U\|_{\sbmod}\|V\|_{\sbmod},
\end{equation}
\begin{equation}  \label{eq:bipara}
    \| \Delta[U^*,V] \|_{\bmop} \le C \log(n+1) \|U\|_{\bmop}\|V\|_{\bmop},
\end{equation}
\begin{equation}  \label{eq:bimult}
    \| \Delta[U^*,V] \|_{\bmol} \le C (\log(n+1))^2 \| U\|_{\bmol}\| V\|_{\bmol},
\end{equation}
\begin{equation}  \label{eq:binorm}
    \| \Delta[U^*,V] \|_{\bmond} \le C (\log(n+1))^2 \|U\|_{\bmond}\|V\|_{\bmond},
\end{equation}
and the dimensional estimates are sharp.
\end{cor}
\proof Only the upper bounds need to be shown.
 For (\ref{eq:bistrong}), use Proposition \ref{carbmoso} to write $\|B\|_{\sbmod}= \sup_{I\in \DD, \|e\|=1} \|\Lambda_B(h_Ie)\|$
 and (\ref{eq:bidentity}) to estimate
 $$ \| \Delta[U^*,V] \|_{\sbmod}\le \sup_{I\in \DD,\|e\|=1} \| \pi_U^* \pi_V h_I e\|
 + \sup_{I\in \DD,\|e\|=1}\|D_{U^*,V}(h_Ie)\|.$$
Now observe that for $e \in \h$, $I \in \DD$, one has
$$\| \pi_U^* \pi_V h_I e\|   \\
               \le \|U\|_{\bmop} \|V\|_{\sbmod} \|e\| \le C' \log(n+1) \|U\|_{\sbmod}
               \|V\|_{\sbmod}\|e\|
$$
by (\ref{eq:paraest}). Since $D_{U^*,V}
h_Ie = \frac{1}{|I|} \sum_{J\subset I} U_J^* V_Je \, h_I$, one
obtains
\begin{multline*}
 \|D_{U^*,V}(h_Ie)\|= \sup_{f \in \h, \|f\|=1} |\langle D_{U^*,V}(h_Ie), h_I f \rangle| \\
=  \sup_{f \in \h, \|f\|=1}
       \frac{1}{|I|} |\sum_{J\subset I}\langle V_Je, U_Jf\rangle |
   \le \|V_e\|_{\bmodd(\T,\HH)}\|U\|_{\sbmod},
\end{multline*}
and the proof of (\ref{eq:bistrong}) if complete.

Using first (\ref{eq:paraest}) and (\ref{eq:strongmult}) and then
Proposition \ref{prop:bisweep}, we obtain (\ref{eq:bipara}).
In a similar way, using first Proposition \ref{prop:bisweep} and then
(\ref{eq:paraest}), (\ref{eq:strongmult}) yields
(\ref{eq:bimult}).

Finally, for (\ref{eq:binorm}) observe first that for any $U,V \in \opf$, $e, f \in \h$, $t \in \T$,
\begin{eqnarray*}
  && |\langle \Delta[U^*,V](t)e, f \rangle|\\
 &=& |\sum_{I \in \DD} \left \langle \frac{\chi_I(t)}{|I|^{1/2}} V_I e, \frac{\chi_I(t)}{|I|^{1/2}}
                                                                               U_I f \right \rangle| \\
     & \le& \left(\sum_{I \in \DD}  \|\frac{\chi_I(t)}{|I|^{1/2}}V_I e\|^2\right)^{1/2}
          \left(\sum_{I \in \DD} \|\frac{\chi_I(t)}{|I|^{1/2}}U_I f\|^2\right)^{1/2} \\
      &=& \langle S_U (t) e,e \rangle^{1/2} \langle S_V (t) f,f \rangle^{1/2}
     \le\|S_U(t)\|^{1/2} \|S_V(t)\|^{1/2}
\end{eqnarray*}
and therefore
\begin{equation} \label{eq:pointest}
    \|\Delta[U^*,V](t)\| \le \|S_U(t)\|^{1/2} \|S_V(t)\|^{1/2} \quad (t \in \T).
\end{equation}
Now consider the $\bmond$ norm of $\Delta[U^*,V]$. For $I \in \DD$,
\begin{eqnarray*}
   && \|P_I \Delta[U^*,V]\|_{L^1(\T, \mat)}\\
   & =& \|P_I \Delta[P_I U^*,P_I V]\|_{L^1(\T, \mat)} \\
 &\le& 2\|\Delta[P_I U^*,P_I V]\|_{L^1(\T, \mat)} \\
  & \le& 2 \| \|S_{P_I U}(\cdot)\|^{1/2} \|S_{P_I V}(\cdot)\|^{1/2}\|_{L^1(\T)} \\
   &\le& 2 \|S_{P_I U}\|^{1/2}_{L^1(\T,\mat)}\|S_{P_I V}\|^{1/2}_{L^1(\T,\mat)}\\
     &\le & 2 (\log(n+1))^2 \|P_I U\|_{L^2(\T,\mat)} \|P_I U\|_{L^2(\T,\mat)}  \\
     & \le & 2 (\log(n+1))^2 |I| \|U\|_{\bmond} \|V\|_{\bmond},
\end{eqnarray*}
where we obtain the third inequality from (\ref{eq:pointest})
and the fourth inequality from the proof of (\ref{eq:sharpnorm}).
This finishes the proof of (\ref{eq:binorm}). \qed


\begin{thebibliography}{99}
\bibitem[B1]{blasco1} O.~Blasco, \emph{Hardy spaces of vector-valued
functions: Duality},
Trans.~Am.~Math.~Soc. {\bf 308} (1988), no.2, 495-507.
\bibitem[B2]{blasco3} O.~Blasco, \emph{Boundary values of functions in
vector-valued Hardy spaces and geometry of Banach
spaces}, J. Funct. Anal. {\bf 78} (1988), 346-364.

\bibitem[B3]{blasco4} O.~Blasco, \emph{Dyadic BMO, paraproduct and
Haar multipliers}, Contemporary Math. (to appear)

\bibitem[B4]{BlascoArg} O.~Blasco, \emph{
Remarks on operator-valued BMO spaces}, Rev. Uni. Mat. Argentina
{\bf 345} (2004), 63-78.

\bibitem[BP1]{BP} O.~Blasco, S.~Pott, {\em Dyadic BMO on the
bidisk},  Rev. Mat. Iberoamericana {\bf 21} 2(2005), 483-510.

\bibitem[BP2]{new} O.~Blasco, S.~Pott, \emph{Operator-valued dyadic BMO spaces}, to appear in J.~Op.~Th.

\bibitem[Bou] {bourgain}J.~Bourgain, \emph{Vector-valued singular integrals
and the $H^1$-BMO
duality}, Probability Theory and Harmonic Analysis, Cleveland, Ohio 1983
 Monographs and Textbooks in Pure and Applied Mathematics {\bf 98}, Dekker,
New York 1986.

\bibitem[ER]{effros} E.~G.~Effros and Z.~J.~Ruan, Operator Spaces,
London Mathematical Society Monographs 23, Oxford University Press, 2000.
\bibitem[G]{G} A. M. Garsia, \emph{Martingale inequalities: Seminar
Notes on recent progress}, Benjamin, Reading, 1973.
\bibitem[GPTV]{gptv2}T.A.~Gillespie, S.~Pott, S.~Treil, A.~Volberg,
\emph{Logarithmic growth for matrix martingale transforms},
  J. London Math. Soc. (2)  64  (2001),  no. 3, 624-636.

\bibitem[JPa]{jp} B.~Jacob and J.R.~Partington, \emph{The Weiss
conjecture on admissibility of observation operators for contraction semigroups},
Integral Equations and Operator Theory 40 (2001),
231-243.
\bibitem[JPaP]{jpp1} B.~Jacob, J.~R.~Partington, S.~Pott,
\emph{Admissible and weakly admissible observation operators for the right
shift semigroup},
  Proc. Edinb. Math. Soc. (2)  {\bf 45}(2002),  no. 2, 353-362.
\bibitem[K]{katz} N. H. Katz, \emph{Matrix valued paraproducts}, J.
Fourier Anal. Appl. {\bf 300} (1997), 913-921
\bibitem[M]{M} Y. Meyer,\emph{Wavelets and operators} Cambridge
Univ. Press, Cambridge, 1992.
\bibitem[NTV]{ntv} F.~Nazarov, S.~Treil, A.~Volberg,
\emph{Counterexample to
the infinite
dimensional Carleson embedding theorem},
 C. R. Acad. Sci. Paris {\bf 325} (1997), 383-389.
\bibitem[NPiTV]{nptv} F.~Nazarov, G.~Pisier, S.~Treil, A.~Volberg,
 \emph{Sharp estimates in vector Carleson imbedding theorem and for vector
paraproducts}, J. Reine Angew. Math. {\bf 542} (2002), 147-171.
\bibitem[Me1]{mei1} T.~Mei, \emph{Operator valued Hardy spaces}, Memoirs of the Am.~Math.~Soc.~2007, vol.~188, no.~881.
\bibitem[Me2]{mei} T.~Mei, \emph{Notes on Matrix Valued
Paraproducts} Indiana Univ. Math. J. {\bf 55} 2,(2006), 747-760.
\bibitem[Per]{per} M.C.~Pereyra, Lecture notes on dyadic harmonic analysis.
Second Summer School in Analysis and Mathematical Physics (Cuernavaca,
2000), 1--60,
Contemp. Math. {\bf 289},
Amer. Math. Soc., Providence, RI, 2001.
\bibitem[Pet]{petermichl} S.~Petermichl,
\emph{Dyadic shifts and a logarithmic estimate for Hankel operators with
matrix symbol},
  C. R. Acad. Sci. Paris Sér. I Math. {\bf 330}  (2000),  no. 6, 455-460.

 \bibitem[Pi]{Pi1} G. Pisier \emph{ Notes on
Banach spaces valued Hp-spaces, non-commutative martingale
inequalities and related questions}. Preliminary Notes, 2000.

\bibitem[PV]{pvpers} G.~Pisier and A.~Volberg, personal communication
\bibitem[PXu]{pxu} G.~Pisier and Q.~Xu, \emph{Non-commutative martingale
inequalities}, Comm. Math. Physics, 189 (1997) 667-698.
\bibitem[PS]{ps} S.~Pott, C.~Sadosky,
\emph{Bounded mean oscillation on the bidisk and Operator BMO},
J.~Funct.~Anal.~{\bf 189}(2002), 475-495.

\bibitem[PSm]{psm} S.~Pott, M. Smith,
\emph{Vector paraproducts and Hankel operators of Schatten class
via $p$-John-Nirenberg theorem}, J.~Funct.~Anal.~{\bf 217}(2004),  no. 1, 38--78.
\bibitem[SW]{SW} E.~M.~Stein and G.~Weiss, \emph{ Introduction to Fourier
Analysis
on Euclidean Spaces }, Princeton Univ. Press, [1971].
\end{thebibliography}
\end{document}